\documentclass[a4paper,10pt,conference]{ieeeconf}
\IEEEoverridecommandlockouts
\usepackage{graphicx}
\usepackage{amsfonts}
\usepackage{nccmath}
\usepackage{amsmath}
\usepackage{amssymb}
\usepackage{arydshln}
\usepackage{fleqn}
\usepackage{here}
\usepackage{macros_e}
\usepackage{mathfonts}
\usepackage{footnote}

\everymath={\displaystyle}
\title{
\LARGE \bf 
$L_{2+}$ Induced Norm Analysis of Continuous-Time LTI Systems\\
Using Positive Filters and Copositive Programming
}
\author{Yoshio Ebihara, Hayato Waki, Noboru Sebe,\\ 
Victor Magron, Dimitri Peaucelle, and Sophie Tarbouriech
\thanks{
Y. Ebihara is with the 
Graduate School of Information Science and 
Electrical Engineering, Kyushu University, 
744 Motooka, Nishi-ku, Fukuoka 819-0395, Japan, 
he was also with 
LAAS-CNRS, Universit\'{e} de Toulouse, CNRS, Toulouse, France,   
in 2011.  
H. Waki is with the Institute of Mathematics for Industry,
Kyushu University, 744 Motooka, Nishi-ku, Fukuoka 819-0395, Japan.  
N. Sebe is with the Department of Intelligent and Control
Systems, Kyushu Institute of Technology, Fukuoka 820-8502, Japan.  
V. Magron, D. Peaucelle, and S. Tarbouriech are with 
LAAS-CNRS, Universit\'{e} de Toulouse, CNRS, Toulouse, France.  
}%
}

\begin{document}
\maketitle
\thispagestyle{empty}
\pagestyle{empty}

\begin{abstract}
This paper is concerned with the analysis of the $L_{2}$ induced norm of 
continuous-time LTI systems where the input signals are restricted to be
nonnegative.  This induced norm is referred to as the $L_{2+}$ induced
norm in this paper.  
It has been shown very recently that the $L_{2+}$ induced norm is particularly useful for the
stability analysis of nonlinear feedback systems constructed from
linear systems and static nonlinearities where the nonlinear elements only provide
nonnegative signals.  
For the upper bound computation of the $L_{2+}$ induced norm, 
an approach with copositive programming has also been proposed.  
It is nonetheless true that this approach becomes effective only for multi-input systems, 
and for single-input systems
this approach does not bring any improvement 
over the trivial upper bound, the standard $L_2$ norm.  
To overcome this difficulty, we newly introduce 
positive filters to increase the number of positive signals.  
This enables us to enlarge the size of the copositive multipliers 
so that we can obtain better (smaller) upper bounds with copositive programming.  

\end{abstract}


\section{Introduction}

Recently, control theoretic approaches for the analysis of 
control systems driven by neural networks (NNs)
have attracted great attention 
\cite{Anderson_IEEENN2007,Yin_IEEE2022,Revay_LCSS2021}.  
The basic treatment there is to recast the control system of interest
into a feedback system constructed from a linear system and
nonlinear activation functions.  
Then, by  capturing the properties of 
nonlinear activation functions with the 
integral quadratic constraint (IQC) framework \cite{Megretski_IEEE1997},   
we can obtain numerically tractable
semidefinite programming problems (SDPs) for the analysis.  

On the other hand, some activation functions in NNs exhibit 
particular nonnegative properties
that are hardly captured by the standard IQC framework on the 
positive semidefinite cone.  
As a typical example, the rectified linear units (ReLUs) return
only nonnegative output signals irrespective of input signals.  
This is the motivation of \cite{Motooka_ISCIE2022}
to consider the $L_2$ induced norm 
of continuous-time LTI systems where the input signals are restricted 
to be nonnegative.  
This induced norm is referred to as the $L_{2+}$ induced
norm in this paper.  
As the main result of \cite{Motooka_ISCIE2022}, 
an $L_{2+}$-induced-norm-based small gain theorem has been derived  
for the stability analysis of recurrent neural networks 
with activation functions being ReLUs.   
This has been proved to be less conservative
than the standard $L_{2}$-induced-norm-based small gain theorem \cite{Khalil_2002}.  
Moreover, for the upper bound computation of the $L_{2+}$ induced norm, 
an approach with copositive programming (COP) has also been proposed in 
\cite{Motooka_ISCIE2022}.  
By applying inner approximation to the copositive cone, 
we can eventually obtain numerically tractable SDPs for the upper bound computation.  

Even though the preceding work \cite{Motooka_ISCIE2022} provides basic ideas 
for the treatment of the $L_{2+}$ induced norm, 
the results there are certainly deficient 
in the following aspects:
\begin{itemize}
 \item[(i)] For single-input systems, the upper bound characterized by the 
	    COP in \cite{Motooka_ISCIE2022} reduces to 
	    the trivial upper bound, the $L_2$ induced norm.  
	    Namely, it is by no means possible to obtain better 
	    upper bounds than the trivial one.  
 \item[(ii)] For multi-input systems, we can obtain a better upper bound than the 
	     $L_2$ induced norm by solving an SDP in \cite{Motooka_ISCIE2022}.  
	     However, once we obtain this upper bound, there is no way 
	     to obtain further better upper bounds.  
\end{itemize}
These deficiencies are related to the size of the
copositive multipliers that is introduced in \cite{Motooka_ISCIE2022}
to capture the nonnegativity of the input signals.  
If we can somehow increase the number of nonnegative signals, 
then we can enlarge the size (freedom) of the corresponding
copositive multiplier so that we can obtain better (smaller) upper bounds.  

To achieve this end, in this paper, we newly introduce 
positive filters to increase the number of nonnegative signals and then 
enlarge the size of copositive multipliers.  
More precisely, we introduce a positive filter of specific form.  
By increasing the degree of the positive filter, 
we can construct a sequence of COPs and 
then a sequence of SDPs by applying inner approximation to the copositive cone.  
We prove that, by solving the sequence of SDPs, 
we can construct a monotonically nonincreasing sequence of the upper bounds 
of the $L_{2+}$ induced norm.  
The effectiveness of the proposed positive-filter-based method
is illustrated by numerical examples.  
We finally note that 
the analysis of the $L_{2+}$ induced norm is also motivated
by recent advancement on the study of positive systems
\cite{Briat_IJRN2013,Tanaka_IEEE2011,Rantzer_IEEE2016,Ebihara_IEEE2017,Kato_LCSS2020}, where the treatment of nonnegative signals is essentially important.   

Notation: 
The set of $n\times m$ real matrices is denoted by $\bbR^{n\times m}$, and
the set of $n\times m$ entrywise nonnegative matrices is denoted
by $\bbR_+^{n\times m}$. 
For a matrix $A$, we also write $A\geq 0\ (A>0)$ to denote that  
$A$ is entrywise nonnegative (positive). 
We denote the set of $n\times n$ real symmetric, 
positive semidefinite, and positive definite matrices by $\bbS^n$, $\bbS_+^n$, and
$\bbS_{++}^n$, respectively.  
The set of $n\times n$ Hurwitz and Metzler matrices are denoted by 
$\bbH^n$ and $\bbM^n$, respectively.  
For $A\in\bbS^n$, we also write $A\succ 0\ (A\prec 0)$ to
denote that $A$ is positive (negative) definite.  

\section{Preliminaries}

\subsection{Norms for Signals and Systems}

For a continuous-time signal $w$ defined 
over the time interval $[0,\infty)$, we define
\[
\begin{array}{@{}l}
 \|w\|_{2}:=\sqrt{\int_{0}^{\infty}|w(t)|_2^2 dt}
\end{array}
\]
where for $v\in\bbR^{n_v}$ we define
$|v|_2:=\sqrt{\sum_{j=1}^{n_v} v_j^2}$.  
We also define
\[
\begin{array}{@{}l}
 L_{2}  :=\left\{w:\ \|w\|_{2}<\infty \right\},\\ 
 L_{2+} :=\left\{w:\ w\in L_{2},\ w(t)\ge 0\ (\forall t\in[0,\infty)\right\}.    
\end{array}
\]
For an operator $H:\ L_{2}\ni w \mapsto z \in L_{2}$, 
we define its (standard) $L_2$ induced norm by
\begin{equation}
\|H\|_{2}:=\sup_{w\in L_2,\ \|w\|_2=1} \ \|z\|_2.  
\label{eq:L2norm}
\end{equation}
We also define
\begin{equation}
\|H\|_{2+}:=\sup_{w\in L_{2+},\ \|w\|_2=1} \ \|z\|_2.  
\label{eq:L2+norm}
\end{equation}
This is a variant of the $L_2$ induced norm introduced in \cite{Motooka_ISCIE2022}
and referred to as the $L_{2+}$ induced norm in this paper.  
We can readily see that $\|H\|_{2+}\le \|H\|_{2}$.  

\subsection{Copositive Programming}

A copositive programming problem (COP) is a convex optimization problem in which 
we minimize a linear
objective function over the 
linear matrix inequality (LMI) constraints on the copositive cone   
\cite{Dur_2010}.
Even though a COP is a convex optimization problem, 
it is hard to solve it numerically in general.  
We summarize the definitions of cones related to the COP and 
its basics in the appendix section,
where the materials there are borrowed in part from
\cite{Ebihara_EJC2021}.

\subsection{Positive Systems}
\label{sub:positive}

In this paper, we introduce positive filters for the (upper bound) computation 
of the $L_{2+}$ induced norm of continuous-time LTI systems.  
The definition of positivity and related results are briefly summarized 
as follows.  
\begin{definition}\cite{Farina_2000}  
An LTI system is called {\it internally} positive if 
its state and output are both nonnegative for 
any nonnegative input and nonnegative initial state.  
\end{definition}
\begin{definition}\cite{Farina_2000}  
An LTI system is called {\it externally} positive if 
its output is nonnegative for 
any nonnegative input under the zero initial state.  
\end{definition}
\begin{proposition}\cite{Farina_2000}  
An LTI system with coefficient matrices
$A\in\bbR^{n\times n}$, $B\in\bbR^{n\times m}$, 
$C\in\bbR^{l\times n}$ and $D\in\bbR^{l\times m}$
is internally positive if and only if 
\[
A\in\bbM^{n},\ 
B\in\bbR_+^{n\times m},\ 
C\in\bbR_+^{l\times n}, 
D\in\bbR_+^{l\times m} .  
\]
\label{pr:positive}
\end{proposition}
%

\section{$L_{2+}$ Induced Norm Analysis}
\label{sec:L2+}

\subsection{Problem Description}
\label{sub:Problem}

Let us consider the LTI system $G$ given by
\begin{equation}
G:\ 
\left\{
\arraycolsep=0.5mm
\begin{array}{cccccccc}
 \dot x(t)&=& A x(t)& + &B w(t), &\ x(0)=0,\\
   z(t)&=& C x(t)& + &D w(t) \\
\end{array}
\right.  
\label{eq:G}
\end{equation}
where
$A\in\bbR^{n\times n}$, 
$B\in\bbR^{n\times \nw}$,
$C\in\bbR^{\nz\times n}$, and
$D\in\bbR^{\nz\times \nw}$.   
We assume that the system $G$ is stable, i.e., 
the matrix $A$ is Hurwitz stable, and the pair $(A,B)$ is controllable.  
It is well known that the $L_2$ induced norm 
$\|G\|_2$ defined by \rec{eq:L2norm} coincides with the
$H_\infty$ norm for stable LTI systems
and plays an essential role in stability analysis of feedback systems.  
In this paper, we are interested in computing 
the $L_2$ induced norm where the input signal $w$ is 
restricted to be nonnegative.  
Namely, we focus on the computation of 
the $L_{2+}$ induced norm 
$\|G\|_{2+}$ defined by
\rec{eq:L2+norm}.  
As noted, it is very clear that $\|G\|_{2+}\le \|G\|_{2}$.  
Here, it is well known that 
$\|G\|_{2+}= \|G\|_{2}$ holds if $G$ is externally positive, 
see, e.g.,  \cite{Tanaka_IEEE2011,Rantzer_IEEE2016}.  

\subsection{Motivating Example: $L_{2+}$-Induced-Norm-Based Small Gain Theorem}
\label{sub:motivation}

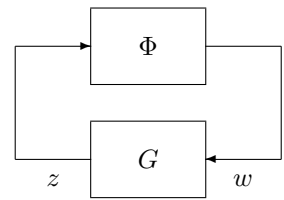
\begin{figure}[h]
\begin{center}
\begin{picture}(3.5,2.5)(0,0)
\put(0,2){\vector(1,0){1}}
\put(1,1.5){\framebox(1.5,1){$\Phi$}}
\put(2.5,2){\line(1,0){1}}
\put(3.5,2){\line(0,-1){1.5}}
\put(3.5,0.5){\vector(-1,0){1}}
\put(3,0.3){\makebox(0,0)[t]{$w$}}
\put(1,0){\framebox(1.5,1){$G$}}
\put(1,0.5){\line(-1,0){1}}
\put(0.5,0.3){\makebox(0,0)[t]{$z$}}
\put(0,0.5){\line(0,1){1.5}}

\end{picture} 
\caption{Nonlinear Feedback System.  }
\label{fig:fb}
\end{center}
\end{figure}

Let us assume $\nz=\nw=m$ for $G$ given by \rec{eq:G}
and consider the feedback system shown in \rfig{fig:fb},   
where $\Phi:\ \bbR^m\mapsto \bbR_+^m$ is a static nonlinear operator 
satisfying $\|\Phi\|_2=1$.  
We focus on the stability analysis of this feedback system.  
Here, note that we have assumed that $\Phi$ returns only nonnegative signals.   
This problem setting typically appears in the stability analysis 
of recurrent neural networks with activation functions being 
rectified linear units, see \cite{Ebihara_EJC2021,Motooka_ISCIE2022}.  

Then, from the standard $L_2$-induced-norm-based 
small-gain theorem \cite{Khalil_2002}, 
we see that the feedback system shown in \rfig{fig:fb} is
(well-posed and) globally stable if $\|G\|_2<1$.  
On the other hand, by actively using the nonnegative nature of $\Phi$, 
it has been shown very recently in \cite{Motooka_ISCIE2022} that 
the feedback system shown in \rfig{fig:fb} is
(well-posed and) globally stable if $\|G\|_{2+}<1$.  
As illustrated by this concrete example, 
the $L_{2+}$-induced-norm-based small-gain theorem
has potential abilities for the stability analysis of feedback systems
with nonnegative nonlinearities.   
This strongly motivates us to establish 
efficient methods for the (upper bound) computation of the 
$L_{2+}$ induced norm.  

\subsection{Basic Results}
\label{sub:basic}

The next result forms an important basis of this study.  
\begin{proposition}\cite{Motooka_ISCIE2022}
For the LTI system $G$ given by \rec{eq:G} and a given $\gamma>0$, 
let us consider the following conditions (i) and (ii).  
\begin{itemize}
 \item[(i)]  $\|G\|_{2+}\le \gamma$.  
 \item[(ii)]  There exist $P\in \PSD^{n}$ and $Q\in \COP^{\nw}$ such that
\begin{equation}
\begin{array}{@{}l}
\left[
\arraycolsep=0.5mm
\begin{array}{cc}
 PA+A^TP & PB+C^TD \\
 \ast & D^TD-\gamma^2 I_{\nw}+Q
\end{array}
\right]\preceq 0.  
\end{array}
\label{eq:LMI}
\end{equation}
\end{itemize}
Then we have (i) $\Leftarrow$ (ii).  
\label{pr:basic}
\end{proposition}

The copositive matrix variable $Q$ in \rec{eq:LMI}
is referred to as the copositive multiplier in this paper.  
On the basis of \rpr{pr:basic}, let us consider the COP:
\begin{equation}
 \ogam:=\inf_{\gamma,P,Q}\ \mbox{subject to}\ \rec{eq:LMI},\ 
  P\in\PSD^n,\ 
  Q\in\COP^{\nw}.  
\label{eq:COP}
\end{equation}
In relation to this COP, recall that
\[
 \|G\|_2=\inf_{\gamma,P}\ \mbox{subject to}\ \rec{eq:LMI},\ 
  P\in\PSD^n,\ 
  Q=0.  
\]
It follows that $\|G\|_{2+}\le \ogam\le \|G\|_2$.  
Unfortunately, as we have already mentioned, 
it is hard to solve the COP \rec{eq:COP} in general.  
However, an upper bound of $\ogam$
can be computed efficiently by replacing the copositive cone $\COP$ in 
\rec{eq:COP} with the Minkowski sum of the
positive semidefinite and nonnegative cones
$\PSD+\NN$ as follows:  
\begin{equation}
\begin{array}{@{}l}
 \oogam:=\inf_{\gamma,P,Q}\ \mbox{subject to}\ \rec{eq:LMI},\\ 
  P\in\PSD^n,\ 
  Q\in\PSD^{\nw}+\NN^{\nw}.  
\end{array}
\label{eq:SDP}
\end{equation}
Note that this problem is essentially an SDP and hence tractable.  
We can readily see that 
$\|G\|_{2+}\le \ogam\le \oogam\le \|G\|_2$ holds.    

Up to this point, we have described the basic ideas of
the (upper bound) computation of $\|G\|_{2+}$ given in \cite{Motooka_ISCIE2022}.  
However, these results are certainly deficient in the following aspects:
\begin{itemize}
 \item[(i)] In the case where $\nw=1$, i.e., 
	    if the system $G$ has only a single disturbance input, 
	    then it is very clear that $\ogam=\|G\|_{2}$.  
	    This is because, since
	    $\COP^1=\PSD^1=\bbR_+$, and since the copositive multiplier $Q$ enters
	    in block-diagonal part in \rec{eq:LMI}, 
	    we see that the optimal value of $Q$ in COP \rec{eq:COP} is zero.  
	    Namely, if $\nw=1$, 
	    it is impossible to obtain
	    an upper bound of $\|G\|_{2+}$ 
	    which is better than the trivial upper bound  
	    $\|G\|_{2}$ if we directly work on \rec{eq:COP}.  
 \item[(ii)] In the case where $n_w>1$, it has been shown by numerical examples 
	     in \cite{Motooka_ISCIE2022}
	     that we can obtain an upper bound $\oogam$
	     of $\|G\|_{2+}$ that is strictly better than 
	     $\|G\|_2$ (i.e., $\oogam<\|G\|_2$).  
	     However, once we obtain $\oogam$, 
	     we have no way to obtain further better upper bounds than $\oogam$.  
\end{itemize}

To overcome these difficulties, we provide a new way for the upper bound computation
of the $L_{2+}$ induced norm with positive filters in this paper.  
\begin{remark}
In \cite{Ebihara_EJC2021},   
the $l_{2+}$ induced norm for 
discrete-time operators has been defined in a similar way to \rec{eq:L2+norm}.   
In the COP-based treatments of 
$l_{2+}$ induced norm computation of discrete-time LTI systems, 
we of course encounter the same difficulties as (i) and (ii) described above.  
However, in the case of discrete-time systems, 
we can employ the discrete-time system lifting \cite{Bittanti_1996}
so that we can artificially increase the number of disturbance inputs.  
This enables us to employ copositive multipliers of larger size.  
By means of this lifting-based treatment, 
we can establish an effective method for the upper bound computation of 
$l_{2+}$ induced norm, see \cite{Ebihara_EJC2021} for details.  
We finally note that there is no genuine counterpart of the lifting in the continuous-time
system setting, and hence we surely need an alternative solution.  
\end{remark}
%

\section{Better Upper Bound Computation by \\Positive Filters}
\label{sec:upper}

For better upper bound computation of $\|G\|_{2+}$, 
it is promising to actively use the fact that the input signal
$w$ is restricted to be nonnegative.    
To this end, let us introduce the positive filter given by
\begin{equation}
G_p : \begin{cases}
 \dot x_p(t)=A_p x_p(t)+B_pw(t),\ x_p(0)=0,\\
 z_p(t)= 
 \left[
  \begin{array}{c}
   I_{\np}  \\
   0_{\nw,\np}
  \end{array}
 \right] x_p(t)+
 \left[
  \begin{array}{c}
   0_{\np,\nw}\\
   I_{\nw}  \\
  \end{array}
 \right] w(t).\hspace*{-20mm}
\end{cases}
\label{eq:Gp}
\end{equation}
Here, $A_p\in \bbH^{\np}\cap \bbM^{\np}$, $B_p\in\bbR_+^{\np\times \nw}$.  
It is clear that the filter $G_p$ is positive from \rpr{pr:positive}.  

We next plug $G_p$ with $G$ and construct the augmented system $G_a$
given by
\begin{equation}
 \scalebox{0.95}{$
 \begin{array}{@{}l}
G_a : \begin{cases}
\dot x_a(t) = A_a x_a(t)+ B_a w(t),\\
z(t)= C_a x_a(t) +D_aw(t),\\
 z_p(t)= 
 \left[
  \begin{array}{cc}
   0_{\np,n} & I_{\np}  \\
   0_{\nw,n} & 0_{\nw,\np}
  \end{array}
 \right]x_a(t)+
 \left[
  \begin{array}{c}
   0_{\np,\nw}\\
   I_{\nw}  \\
  \end{array}
 \right] w(t).  \hspace*{-20mm}
\end{cases}
 \end{array}$}
\label{eq:Ga}
\end{equation}
Here, 
\begin{equation}
\begin{array}{@{}l}
x_a:=
\left[
  \begin{array}{c}
   x \\
   x_p \\
  \end{array}
\right],\ 
A_a:= 
\left[
  \begin{array}{cc}
   A & 0  \\
   0 & A_p \\
  \end{array}
\right],\ 
B_a:=
\left[
  \begin{array}{c}
   B \\
   B_p \\
  \end{array}
\right],\\
C_a:=
\left[
  \begin{array}{cc}
   C & 0_{n_z,\np}  \\
  \end{array}
\right], D_a:=D.  
\end{array}
\label{eq:aug}
\end{equation}
In the above augmented system $G_a$, 
it is very important to note that 
the output $z_p$ is nonnegative for any nonnegative 
input $w$.  
By focusing on this property, we can obtain 
the first main result of this paper as summarized in the next theorem.  
\begin{theorem}
For the LTI system $G$ given by \rec{eq:G} and a given $\gamma>0$, 
let us consider the following conditions (i) and (iii).  
\begin{itemize}
 \item[(i)]  $\|G\|_{2+}\le \gamma$.  
 \item[(iii)]  There exist $P_a \in \bbS^{n+\np}$ and 
	       $Q_a\in \COP^{\np+\nw}$ such that
\begin{equation}
\begin{array}{@{}l}
\begin{bmatrix}
P_aA_a + A_a^T P_a + C_a^T C_a & P_aB_a + C_a^T D_a \\
B_a^T P_a + D_a^T C_a & D_a^T D_a- \gamma ^2 I_{\nw}
\end{bmatrix}\\
+
\begin{bmatrix}
0_{n,\np+\nw} \\ I_{\np+\nw}
\end{bmatrix}
Q_a
\begin{bmatrix}
0_{n,\np+\nw} \\ I_{\np+\nw}
\end{bmatrix}^T
\preceq 0.  
\end{array}
\label{eq:L2+COPnew}
\end{equation}
\end{itemize}
Then, we have (i) $\Leftarrow$ (iii).  
\label{th:main1}
\end{theorem}

For the proof of this theorem, we need the next lemmas.  
In the following, we partition $Q_a\in\COP^{\np+\nw}$ as
\[
\arraycolsep=0.5mm
 Q_{a}=
 \left[
 \begin{array}{cc}
  Q_{a,11} & Q_{a,12}\\
  Q_{a,12}^T & Q_{a,22}\\
 \end{array}
 \right],\ Q_{a,11}\in\COP^{\np},\ Q_{a,22}\in\COP^{\nw}.  
\]
\begin{lemma}
For $A_p\in\bbH^{\np}\cap\bbM^{\np}$ and $Q_{a,11}\in\COP^{\np}$, 
let us consider the unique solution $P_p\in\bbS^{\np}$ to the Lyapunov equation
\begin{equation}
 P_p A_p + A_p^T P_p +Q_{a,11}=0.   
\label{eq:Lyap_Ap}
\end{equation}
Then, we have $P_p\in\COP^{\np}$.  
\label{le:pos_Lyap}
\end{lemma}
\begin{proofof}{\rle{le:pos_Lyap}}
The proof can readily be done if we note that
\[
 P_p=\int_0^\infty \exp(A_p^Tt)Q_{a,11}\exp(A_pt)dt
\]
and $\exp(A_pt)\ge 0\ (\forall t\ge 0)$ holds since $A_p$ is Metzler.  
\end{proofof}
\begin{lemma}
Suppose $P_a\in\bbS^{n+\np}$ satisfies \rec{eq:L2+COPnew} with
$Q_a\in \COP^{\np+\nw}$.  Then, we have
\begin{equation}
\arraycolsep=0.5mm
  \scalebox{0.87}{$
  \begin{array}{@{}l}
   P_a\in \left\{P+\begin{bmatrix}
		    0_{n,n} & 0_{n,\np} \\
		    0_{\np,n} & P_p \\
		   \end{bmatrix}:\ P\in \PSD^{n+\np},\ P_p\in \COP^{\np}\right\}.  
  \end{array}$}
\label{eq:Pa}
\end{equation}
\label{le:structure}
\end{lemma}
\begin{proofof}{\rle{le:structure}}
Suppose $P_a\in\bbS^{n+\np}$ satisfies \rec{eq:L2+COPnew} with
$Q_a\in \COP^{\np+\nw}$.  
Then, it is very clear that there exists $W\in\bbS_+^{n+\np}$ such that
\[
 P_a A_a + A_a^T P_a +C_a^T C_a + 
 \left[
 \begin{array}{cc}
  0_{n,n} & 0_{n,\np} \\
  0_{\np,n} & Q_{a,11}
 \end{array}
 \right]+W=0   
\]
where $A_a\in\bbH^{n+\np}$ from \rec{eq:aug}.  
If we regard this equation as the Lyapunov equation 
with respect to $P_a\in\bbS^{n+\np}$, we see from the linearity that
$P_a\in\bbS^{n+\np}$ can be written as 
\[
 P_a=P+\begin{bmatrix}
		    0_{n,n} & 0_{n,\np} \\
		    0_{\np,n} & P_p \\
		   \end{bmatrix}.  
\]
Here, $P\in\bbS_+^{n+\np}$ is the unique solution to the Lyapunov equation
\[
 P A_a + A_a^T P +C_a^T C_a + W=0 
\]
whereas $P_p\in\bbS^{\np}$ is the unique solution to the Lyapunov equation
\rec{eq:Lyap_Ap}.   From \rle{le:pos_Lyap}, 
we have $P_p\in\COP^{\np}$ and hence   
the proof is completed.  
\end{proofof}

Now we are ready to prove \rth{th:main1}.  

\begin{proofof}{\rth{th:main1}}
For the augmented system $G_a$, we consider the trajectory of its state $x_a$
for the input $w\in L_{2+}$ with $\| w \| _2 = 1$.  
From \rec{eq:L2+COPnew}, we readily see
\[
\scalebox{0.95}{$
\begin{array}{@{}l}
\begin{bmatrix}
x_a(t) \\ w(t)
\end{bmatrix}
^\mathrm{T}
\begin{bmatrix}
P_aA_a + A_a^T P_a + C_a^T C_a & P_aB_a + C_a^T D_a \\
B_a^T P_a + D_a^T C_a & D_a^T D_a- \gamma^2 I_{\nw}
\end{bmatrix}
\begin{bmatrix}
x_a(t) \\ w(t)
\end{bmatrix}\\
+
\begin{bmatrix}
x_p(t) \\ w(t)
\end{bmatrix}^T
Q_a
\begin{bmatrix}
x_p(t) \\ w(t)
\end{bmatrix}
\le 0\ (\forall t \ge 0).  
\end{array}$}
\]
From this inequality and \rec{eq:Ga}, we have
\[
\begin{array}{@{}l}
\frac{d}{dt} x_a(t)^T P_a x_a(t) + z(t)^T z(t) - \gamma^2  w(t)^T w(t) \\
 +
z_p(t)^T Q_a z_p(t)
\le 0\ (\forall t \ge 0).  
\end{array}
\]
By integration, we arrive at
\begin{equation}
\begin{array}{@{}l}
 x_a(T)^T P_a x_a(T) + \int_{0}^T z(t)^T z(t) - \gamma^2  w(t)^T w(t) dt\\
 +
\int_{0}^T z_p(t)^T Q_a z_p(t) dt
\le 0\ (\forall T > 0).  
\end{array}
\label{eq:int}
\end{equation}
Since $Q_a \in \COP^{\np+\nw}$ and $z_p\in L_{2+}$, we first note that
\begin{equation}
\int_{0}^T z_p(t)^T Q_a z_p(t) dt\ge 0\ (\forall T>0).   
\label{eq:hard_IQC}
\end{equation}
On the other hand, 
since $x_p$ is nonnegative in $x_a(=[\ x^T\ x_p^T\ ]^T)$, 
we see from \rle{le:structure} that 
\begin{equation}
x_a(T)^T P_a x_a(T) \ge 0\ (\forall T>0) .  
\label{eq:quadpos}
\end{equation}
With these facts in mind, we take the limit $T\to\infty$ in \rec{eq:int}
and obtain
\begin{equation}
 \int_0^{\infty}  z(t)^T z(t) \:dt - \gamma^2 \int_0^{\infty} w(t)^T w(t) \:dt \le 0.  
\label{eq:key} 
\end{equation}
This clearly shows $ \| z \|_2 ^2 \le \gamma^2 \| w \|_2 ^2 = \gamma^2$.  
To summarize, we arrive at the conclusion that
\[
\| G \| _{2+} = \sup_{w \in L_{2+}, \| w \| _2 = 1} \| z \| _2 \ \le \ \gamma.  
\]
This completes the proof.  
\end{proofof}
\begin{remark}
In stark contrast with the standard $L_2$-induced norm computation case, 
the Lyapunov certificate 
$P_a \in \bbS^{n+\np}$ that satisfies \rec{eq:L2+COPnew} does not satisfy
$P_a \in \PSD^{n+\np}$ in general.  Namely, we have \rec{eq:Pa}.  
\label{re:nonpos}
\end{remark}
\begin{remark}
Even though our proof of \rth{th:main1} 
has been done by purely time-domain arguments, 
it has close relationship with the 
general IQC theorem \cite{Scherer_Automatica2008,Scherer_SCL2018}
basically characterized in frequency domain.    
Recall that the key conditions in the proof are \rec{eq:hard_IQC} and 
\rec{eq:quadpos}, and these are related to the motivation of 
the study in \cite{Scherer_SCL2018}.    
The key condition \rec{eq:hard_IQC} is referred to 
as the hard (finite horizon) IQC in \cite{Scherer_SCL2018}, 
and its weaker condition is called the soft (infinite horizon) IQC.  
Moreover, even though \rec{eq:quadpos} does hold, 
the Lyapunov certificate $P_a$ in \rec{eq:L2+COPnew} 
is not positive semidefinite in general (\rle{le:structure}),  
and this conforms to the general setting of \cite{Scherer_SCL2018}
that deals with indefinite Lyapunov certificates and soft  (infinite horizon) IQCs.  
Therefore, it seems that the current result can be 
basically subsumed into the general framework of \cite{Scherer_SCL2018}.  
However, since we have introduced COP multipliers, 
the relationship of \rth{th:main1}
with \cite{Scherer_SCL2018} is not yet very clear to us.   
This topic is currently under investigation.  
\end{remark}

On the basis of \rth{th:main1}, let us consider the following 
COP and SDP:  
\begin{equation}
 \scalebox{1.0}{$
\begin{array}{@{}l}
\ogam_a := \inf_{\gamma, P_a, Q_a}\ \gamma\quad \mathrm{subject\ to}\  
\rec{eq:L2+COPnew},\\ 
P_a \in \bbS^{n+\np},\ Q_a \in \COP^{\np+\nw}, 
\end{array}$}
\label{eq:COPnew}
\end{equation}
\begin{equation}
 \scalebox{1.0}{$
\begin{array}{@{}l}
\oogam_a := \inf_{\gamma, P_a, Q_a}\ \gamma\quad \mathrm{subject\ to}\
 \rec{eq:L2+COPnew},\\ 
P_a \in \bbS^{n+\np},\ Q_a \in \PSD^{\np+\nw}+\NN^{\np+\nw}.  
\end{array}$}
\label{eq:SDPnew}
\end{equation}
Then, we readily obtain 
\begin{equation}
\|G\|_{2+} \le \ogam_a \le \oogam_a.   
\end{equation}
Moreover, we can obtain the next theorem 
that verifies the effectiveness of the introduction of positive filters 
in computing better (smaller) upper bounds.  
\begin{theorem}
Let us consider the positive-filter-based 
upper bounds $\ogam_a$ and $\oogam_a$
of $\|G\|_{2+}$ given respectively   
by \rec{eq:COPnew} and \rec{eq:SDPnew}.    
Then, in relation to the filter-free upper bounds $\ogam$ and $\oogam$
given respectively by \rec{eq:COP} and \rec{eq:SDP}, 
we have
\begin{equation}
\ogam_a\le \ogam\ (\le \|G\|_2),\quad \oogam_a\le \oogam\  (\le \|G\|_2)
\label{eq:impr}
\end{equation}
\label{th:impr}
\end{theorem}
\begin{proofof}{\rth{th:impr}}
In the following, we prove $\ogam_a\le \ogam$.  
The proof for $\oogam_a\le \oogam$ follows similarly.  
To prove $\ogam_a\le \ogam$, 
it suffices to show that   
the condition \rec{eq:L2+COPnew} in (iii) of
\rth{th:main1} holds with $\gamma=\ogam+\varepsilon$
for any $\varepsilon>0$.  

We first note from the definition of $\ogam$
given by \rec{eq:COP} that 
for any $\varepsilon>0$
there exist 
$P\in\bbS^n$,  $Q\in\COP^{\nw}$, and $\varepsilon_1>0$ such that
\[
\scalebox{0.77}{$
\begin{array}{@{}l}
\begin{bmatrix}
PA + A^T P + C^T C & PB + C^T D & 0 \\
\ast & D^T D- (\ogam+\varepsilon)^2 I_{\nw}+Q & \varepsilon_1 B_p^T P_p\\
\ast & \ast & \varepsilon_1 (P_pA_p+A_p^TP_p)
\end{bmatrix}
\preceq 0.  
\end{array}$}
\]
Here, $P_p\in\bbS_{++}^{\np}$ is the unique solution of the Lyapunov equation
 \begin{equation}
 P_pA_p+A_p^TP_p+I=0.  
\label{eq:Lyapp}
 \end{equation}
By applying a congruence transformation to the preceding inequality, we have
\[
\scalebox{0.77}{$
\begin{array}{@{}l}
\begin{bmatrix}
PA + A^T P + C^T C & 0 & PB + C^T D  \\
\ast & \varepsilon_1 (P_pA_p+A_p^TP_p) & \varepsilon_1 P_p B_p \\
\ast & \ast & D^T D- (\ogam+\varepsilon)^2 I_{\nw}+Q 
\end{bmatrix}
\preceq 0.  
\end{array}$}
\]
This clearly shows that \rec{eq:L2+COPnew} in (iii) of \rth{th:main1}
holds with $\gamma=\ogam+\varepsilon$ and
\[
\scalebox{0.9}{$
\begin{array}{@{}l}
 P_a=
 \begin{bmatrix}
  P & 0 \\ 0 & \varepsilon_1 P_p
  \end{bmatrix}\in\bbS^{n+\np},\ 
 Q_a=
 \begin{bmatrix}
  0_{\np,\np} & 0 \\ 0 & Q
  \end{bmatrix}\in\COP^{\np+\nw}.  
\end{array}$}
\]
This completes the proof.  
\end{proofof}

\rth{th:impr} shows that 
for any positive filter $G_p$ the corresponding upper bound
$\ogam_a$ given by \rec{eq:COPnew} ($\oogam_a$ given by \rec{eq:SDPnew})
is better (no worse) than the filter-free upper bound 
$\ogam$ given by \rec{eq:COP} ($\oogam$ given by \rec{eq:SDP}).  
In \rsec{sec:num}, we demonstrate by numerical examples that
$\oogam_a$ can be strictly better than $\oogam$ as expected.  

\section{Concrete Construction of Positive Filters}
\label{sec:filter}

As for the positive filter $G_p$ given by \rec{eq:Gp}, 
let us consider the specific form given by
\begin{equation}
\scalebox{0.9}{$
\begin{array}{@{}l}
A_p=A_{p,\alpha,N}:=J_{\alpha,N}\otimes I_{\nw}\in\bbR^{N\nw\times N\nw},\\
B_p=B_{p,N}:=E_N\otimes I_{\nw}\in\bbR^{N\nw\times\nw}, \\
J_{\alpha,N}:=
 \left[
 \begin{array}{ccccc}
  \alpha & 1 & 0 & \cdots & 0 \\
  0 & \alpha & 1 &  \ddots & \vdots \\
  \vdots & \ddots & \ddots &  \ddots & 0 \\
  \vdots &   & \ddots &  \ddots & 1 \\
  0 & \cdots & \cdots &  0 & \alpha \\
 \end{array}
 \right]\in \bbR^{N\times N},\
E_N:=
 \left[
 \begin{array}{c}
  0 \\
  \vdots\\
  0 \\
  1
 \end{array}
 \right]\in \bbR^{N} \hspace*{-20mm} 
\end{array}$}
\label{eq:FIR}
\end{equation}
where $\alpha<0$.  
In this case, the input-output property of
the positive filter $G_p$ in frequency domain is given by
\[
 Z_p(s)=
\left[
\begin{array}{c}
 \dfrac{1}{(s-\alpha)^N}I_{\nw}\\
 \vdots\\
 \dfrac{1}{(s-\alpha)}I_{\nw}\\
 I_{\nw}
\end{array}
\right]W(s).  
\]

By increasing the degree $N$ of the positive filter $G_p$ given 
by \rec{eq:Gp} and \rec{eq:FIR}, 
we can construct a sequence of COPs in the form of \rec{eq:COPnew}
and SDPs in the form of \rec{eq:SDPnew}.  
In the following, we denote by 
$\ogam_{a,\alpha,N}$ and $\oogam_{a,\alpha,N}$
the optimal values of these COPs and SDPs, respectively.  
In addition, we denote by 
$A_{a,\alpha,N}$, 
$B_{a,N}$, 
$C_{a,N}$, and 
$D_{a,N}(=D)$
the coefficient matrices of the augmented system $G_a$
given by \rec{eq:Ga} corresponding to the filter of degree $N$.  
Then, regarding the effectiveness of employing higher-degree positive filters
in improving upper bounds, 
we can obtain the next result.  
\begin{theorem}
Let us consider the upper bounds of $\|G\|_{2+}$
given by $\ogam_{a,\alpha,N}$ and $\oogam_{a,\alpha,N}$ 
that are characterized respectively by 
\rec{eq:COPnew} and \rec{eq:SDPnew} with 
the positive filter $G_p$ of the form 
\rec{eq:Gp} and \rec{eq:FIR} of degree $N$.  
Then, for $N_1\le N_2$, we have
\begin{equation}
\ogam_{a,\alpha,N_2}\le \ogam_{a,\alpha,N_1},\quad 
\oogam_{a,\alpha,N_2}\le \oogam_{a,\alpha,N_1}.    
 \label{eq:monotone}
\end{equation}
 \label{th:monotone}
\end{theorem}
\begin{proofof}{\rth{th:monotone}}
In the following, we prove $\ogam_{a,\alpha,N_2}\le \ogam_{a,\alpha,N_1}$.  
The proof for $\oogam_{a,\alpha,N_2}\le \oogam_{a,\alpha,N_1}$
follows similarly.  
To prove $\ogam_{a,\alpha,N_2}\le \ogam_{a,\alpha,N_1}$, 
it suffices to show that 
$\ogam_{a,\alpha,N+1}\le \ogam_{a,\alpha,N}$ holds for any $N$.  
Furthermore, this can be verified by proving that 
\rec{eq:L2+COPnew} corresponding to the filter of degree $N+1$
holds with $\gamma=\ogam_{a,\alpha,N}+\varepsilon$ for any 
$\varepsilon>0$.  

To this end, we first note from the definition of 
$\ogam_{a,\alpha,N}$ that 
for any $\varepsilon>0$ there exist 
$P_a=P_{a,\alpha,N}\in \bbS^{n+N\nw}$ and  
$Q_a=Q_{a,\alpha,N}\in \COP^{(N+1)\nw}$ such that 
\begin{equation}
 \scalebox{0.95}{$
\begin{array}{@{}l}
\begin{bmatrix}
P_{a,\alpha,N}A_{a,\alpha,N} + A_{a,\alpha,N}^T P_{a,\alpha,N} & P_{a,\alpha,N}B_{a,N} \\
B_{a,N}^T P_{a,\alpha,N}  & - (\ogam_{a,\alpha,N}^2+2\ogam_{a,\alpha,N}\varepsilon) I_{\nw}
\end{bmatrix}\hspace*{-20mm}\\
+
\begin{bmatrix}
C_{a,N}^T \\ D_{a,N}^T
\end{bmatrix}
\begin{bmatrix}
C_{a,N}^T \\ D_{a,N}^T
\end{bmatrix}^T\\
+
\begin{bmatrix}
0_{n,(N+1)\nw} \\ I_{(N+1)\nw}
\end{bmatrix}
Q_{a,\alpha,N}
\begin{bmatrix}
0_{n,(N+1)\nw} \\ I_{(N+1)\nw}
\end{bmatrix}^T
\preceq 0.  
\end{array}$}
\label{eq:N}
\end{equation}
To proceed, let us define
\[
\begin{array}{@{}l}
F_N:=
\begin{bmatrix}
 I_{\nw} & 0_{\nw,\nw} & \cdots & 0_{\nw,\nw}
 \end{bmatrix}\in\bbR^{\nw\times N\nw}.  
\end{array}
\]
Then, there exist $\varepsilon_1,\varepsilon_2>0$ such that
\[
\left[
  \begin{array}{ccc}
   0_{n,n}& 0 & 0 \\
   \ast & -\varepsilon_2 I_{N\nw}-\frac{\varepsilon_1}{2\alpha} F_N^TF_N & \varepsilon_2 P_pB_p\\
    \ast & \ast & -\varepsilon^2 I_{\nw}
  \end{array}
\right]\preceq 0  
\]  
where $P_p\in\bbS_{++}^{\np}$ is the unique solution of the 
Lyapunov equation \rec{eq:Lyapp}. 
By summing up the above inequality with \rec{eq:N} and 
applying the Schur complement argument,  
we obtain
\begin{equation}
\scalebox{0.95}{$
\begin{array}{@{}l}
\left[
\begin{array}{cc}
P_{a,\alpha,N}^{11}A + A^T P_{a,\alpha,N}^{11} & 0 \\
 \ast &  2\varepsilon_1 \alpha I_{\nw} \\
\ast & \ast \\
\ast & \ast \\
\end{array}\right.\hspace*{-20mm} \vspace*{2mm}\\
\left.
\begin{array}{cc}
P_{a,\alpha,N}^{12}A_{p,\alpha,N} + A^T P_{a,\alpha,N}^{12} &
P_{a,\alpha,N}^{11}B+P_{a,\alpha,N}^{12}B_p \\
 \varepsilon_1 F_N & 0 \\
 \hatP_{a,\alpha,N}^{22}A_{p,\alpha,N} + A_{p,\alpha,N}^T \hatP_{a,\alpha,N}^{22} & P_{a,\alpha,N}^{12 T}B +\hatP_{a,\alpha,N}^{22}B_p \\
 \ast & -(\ogam_{a,\alpha,N}+\varepsilon)^2 I_{\nw}
\end{array}\right]\hspace*{-20mm} \vspace*{2mm}\\
+
\begin{bmatrix}
C_{a,N+1}^T \\ D_{a,N+1}^T
\end{bmatrix}
\begin{bmatrix}
C_{a,N+1}^T \\ D_{a,N+1}^T
\end{bmatrix}^T\\
+
\begin{bmatrix}
0_{n,(N+2)\nw} \\ I_{(N+2)\nw}
\end{bmatrix}
\begin{bmatrix}
0_{\nw,\nw} & 0 \\ 0 & Q_{a,\alpha,N}
\end{bmatrix}
\begin{bmatrix}
0_{n,(N+2)\nw} \\ I_{(N+2)\nw}
\end{bmatrix}^T
\preceq 0.  
\end{array}$}
\label{eq:N+1}
\end{equation}
Here, 
\[
\scalebox{0.95}{$
\begin{array}{@{}l}
\begin{bmatrix}
P_{a,\alpha,N}^{11} & P_{a,\alpha,N}^{12} \\
P_{a,\alpha,N}^{12 T} & P_{a,\alpha,N}^{22}
\end{bmatrix}:=P_{a,\alpha,N},\
P_{a,\alpha,N}^{11}\in\bbS^{n},\
P_{a,\alpha,N}^{22}\in\bbS^{N\nw},\vspace*{2mm}\\
\hatP_{a,\alpha,N}^{22}:= P_{a,\alpha,N}^{22}+\varepsilon_{2}P_p.  
\end{array}$}  
\]  
Then, \rec{eq:N+1} shows that 
\rec{eq:L2+COPnew} corresponding to the filter of degree $N+1$
holds with $\gamma=\ogam_{a,\alpha,N}+\varepsilon$ and 
\[
 \scalebox{0.95}{$
 \begin{array}{@{}l}
  P_a=P_{a,\alpha,N+1}=
   \begin{bmatrix}
    P_{a,\alpha,N}^{11} & 0 & P_{a,\alpha,N}^{12} \\
    0 & \varepsilon_1 I_{\nw} & 0 \\
    P_{a,\alpha,N}^{12 T} & 0 & \hatP_{a,\alpha,N}^{22}
   \end{bmatrix}\in \bbS^{n+(N+1)\nw},\\   
   Q_a=Q_{a,\alpha,N+1}=
   \begin{bmatrix}
    0_{\nw,\nw} & 0 \\ 0 & Q_{a,\alpha,N}
   \end{bmatrix}\in \COP^{(N+2)\nw}.  
  \end{array}$}  
\]  
This completes the proof.  
\end{proofof}
\begin{remark}
From \rth{th:monotone}, we see that 
we can construct a monotonically non-increasing sequence of 
upper bounds $\{\oogam_{a,\alpha,N}\}$ of $\|G\|_{2+}$
by increasing the degree $N$ and solving the corresponding SDP 
\rec{eq:SDPnew}.  
Of course this is done at the expense of increased computational burden.  
\end{remark}
\begin{remark}
As shown in the appendix section, 
we can prove that 
the (primal) SDP \rec{eq:SDPnew} and its dual both have
interior point solutions (in the case where we employ the positive filter $G_p$
of the form \rec{eq:Gp} and \rec{eq:FIR}).  
Therefore, there is no duality gap between
the SDP \rec{eq:SDPnew} and its dual, and both have optimal solutions, 
see \cite{Klerk_2002}.  
 \end{remark}
%

\section{Numerical Examples}
\label{sec:num}

In this section, we illustrate the effectiveness of the proposed
positive-filter-based method by numerical examples
on single- and multi-input systems.  
We use MOSEK \cite{MOSEK} to solve the SDP \rec{eq:SDPnew}.  
As for the algorithm implemented in MOSEK, 
we can expect the execution of 
reliable computation if the SDP to be solved and its dual 
both have interior point solutions.  
Since the SDP \rec{eq:SDPnew} and its dual indeed both have 
interior solutions, we can conclude that  
we have executed reliable numerical computation in this section.  

\subsection{Single-Input Case}
\label{sub:SI}

\begin{figure}[t]
\vspace*{-3mm}
\begin{center}
\hspace*{-3mm}
 \includegraphics[scale=0.65]{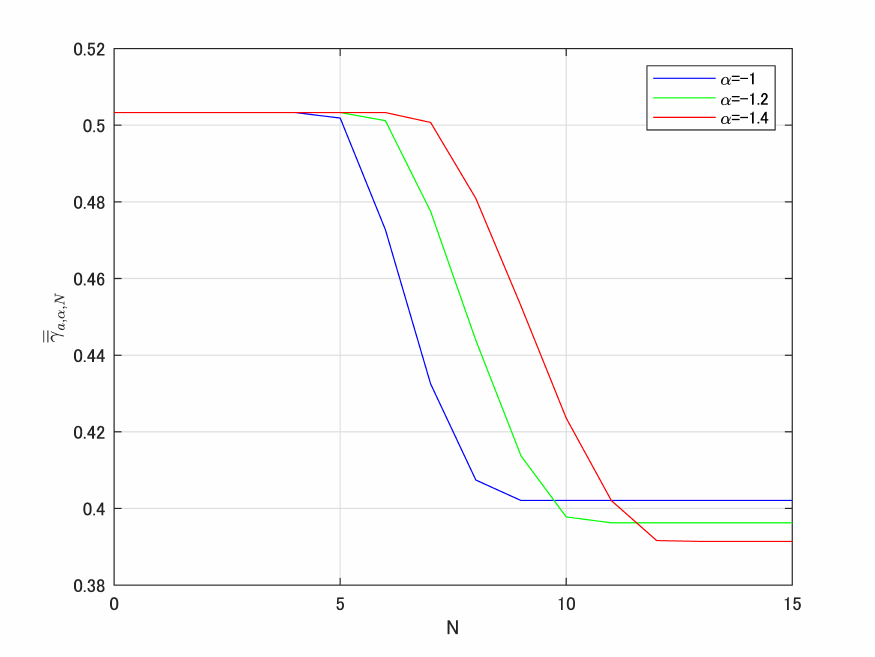}\vspace*{-3mm}\\
 \scalebox{0.8}{The Degree of Positive Filters}\vspace*{-2mm}\\
 \caption{The Values of $\oogam_{a,\alpha,N}$: Upper Bounds of $\|G\|_{2+}$.}
 \label{fig:SI}
\end{center}
\vspace*{-10mm}
\end{figure}

Let us consider the case where the coefficient matrices of the system
\rec{eq:G} are randomly generated and given by
\[
\scalebox{0.95}{$
\begin{array}{@{}l}
 A=
 \begin{bmatrix}
   -0.09 &  0.28 &  0.46 & -0.48 & -0.05\\
  -0.34 & -0.95 & -0.42 &  0.37 & -0.55\\
  -0.24 &  0.04 & -0.10 & -0.47 & -0.23\\
   0.30 &  0.29 &  0.02 & -1.59 &  0.57\\
   0.26 &  0.25 &  0.40 & -0.74 & -0.95\\
 \end{bmatrix},\ 
 B=
 \begin{bmatrix}
 0.17\\
 0.40\\
 0.49\\
 0.30\\
 -0.69\\
 \end{bmatrix},\vspace*{1mm}\\ 
 C=
  \begin{bmatrix}
   -0.14 & -0.66 &  0.10 &  0.34 &  0.05
  \end{bmatrix},\ 
  D=0.27.  
\end{array}$}
\]
In this case, it turned out that $\|G\|_2=0.5033\ (=\oogam_{a,\alpha,0})$.  
Then, for $\alpha\in\{-1,-1.2,-1.4\}$, we constructed 
positive filters of degree $N$ and then 
obtained upper bounds $\oogam_{a,\alpha,N}$  
by solving the SDP \rec{eq:SDPnew}.  
The results are shown in \rfig{fig:SI}.  
We note that $\oogam_{a,\alpha,N}=\ogam_{a,\alpha,N}$ 
holds up to $N=3$ in this case, since for the SDP \rec{eq:SDPnew} we have
$Q_a\in \PSD^{N+1}+\NN^{N+1}=\COP^{N+1}\ (N\le 3)$.  

From \rfig{fig:SI},   
the best (least) upper bound turned out to be 
$\oogam_{a,\alpha,N}=0.3914$ 
with $\alpha=-1.4$ and $N=15$.  
For each $\alpha$, we see that $\oogam_{a,\alpha,N}$
is monotonically non-increasing with respect to $N$ and this result
is surely consistent with \rth{th:monotone}.  

%

\subsection{Multi-Input Case}
\label{sub:MI}

We next consider the case where the coefficient matrices of the system
\rec{eq:G} are randomly generated and given by
\[
\scalebox{0.85}{$
\begin{array}{@{}l}
 A=
 \begin{bmatrix}
  -0.11 & -0.15 &  0.18 &  0.15 & -0.10\\
 0.18 & -0.53 & -0.35 &  0.37 & -0.23\\
  -0.64 & -0.12 & -0.75 &  0.23 &  0.59\\
 0.34 & -0.03 &  0.13 & -0.47 & -0.67\\
 0.55 &  0.29 & -0.08 &  0.53 & -0.81\\
 \end{bmatrix},\ 
 B=
 \begin{bmatrix}
  -0.14 &  0.32\\
  -0.76 & -0.42\\
  -0.30 & -0.03\\
  0.64 & -0.38\\
  -0.12 &  0.17\\
 \end{bmatrix},\vspace*{1mm}\\ 
 C=
  \begin{bmatrix}
   -0.35 &  0.03 &  0.33 &  0.05 &  0.14
  \end{bmatrix},\ 
  D=
  \begin{bmatrix}
 0.43 &  0.23
  \end{bmatrix}.  
\end{array}$}
\]
In this case, it turned out that $\|G\|_2=0.6995$.  
On the other hand, by solving the SDP \rec{eq:SDP} shown in \cite{Motooka_ISCIE2022}, 
we obtained $\oogam=0.6611\ (=\oogam_{a,\alpha,0})$.  
We note that $\oogam=\ogam$ holds in this case 
since for the SDP \rec{eq:SDP} we have
$Q\in \PSD^{2}+\NN^{2}=\COP^{2}$.  

\begin{figure}[b]
\vspace*{-8mm}
\begin{center}
\hspace*{-3mm}
 \includegraphics[scale=0.65]{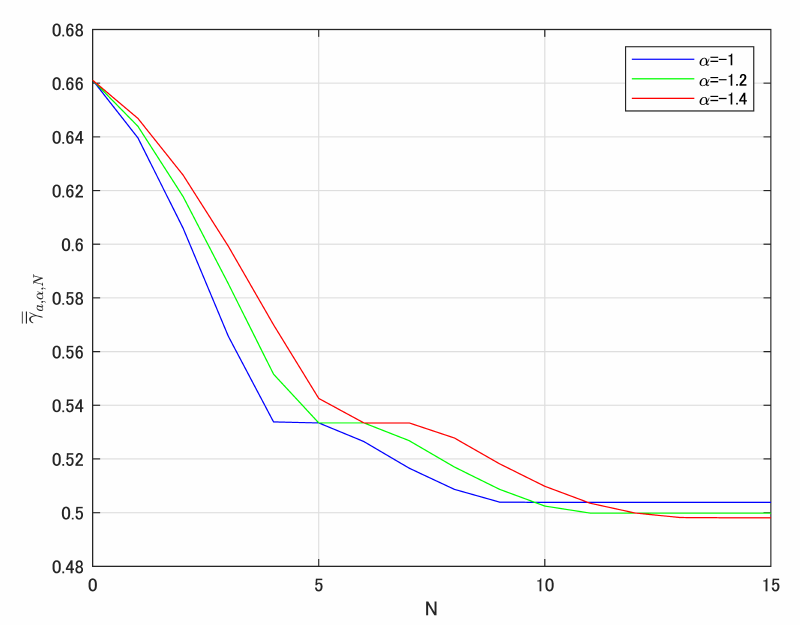}\vspace*{-3mm}\\
 \scalebox{0.8}{The Degree of Positive Filters}\vspace*{-2mm}\\
 \caption{The Values of $\oogam_{a,\alpha,N}$: Upper Bounds of $\|G\|_{2+}$.}
 \label{fig:MI}
\end{center}
\vspace*{-5mm}
\end{figure}

Then, for $\alpha\in\{-1,-1.2,-1.4\}$, we constructed 
positive filters of degree $N$ and then 
obtained upper bounds $\oogam_{a,\alpha,N}$  
by solving the SDP \rec{eq:SDPnew}.  
The results are shown in \rfig{fig:MI}.  
The best (least) upper bound turned out to be 
$\oogam_{a,\alpha,N}=0.4981$ 
with $\alpha=-1.4$ and $N=15$.  
For each $\alpha$, we see that $\oogam_{a,\alpha,N}$
is monotonically non-increasing with respect to $N$ and this result
is again consistent with \rth{th:monotone}.  

%

\section{Conclusion and Future Works}

In this paper, we considered the (upper bond) computation
of the $L_{2+}$ induced norm for continuous-time LTI systems.  
To obtain better (smaller) upper bounds, we introduced 
positive filters and reduced the upper bound computation problem into a COP. 
Then, by applying inner approximation to the COP cone, 
we derived a numerically tractable SDP for the upper bound computation.  
By numerical examples, we showed the effectiveness of 
the proposed positive-filter-based method to obtain  
better (smaller) upper bounds.  

In this paper we just focused on the upper bound computation of the
$L_{2+}$ induced norm.  
It is nonetheless true that we cannot say anything 
on the conservatism of the obtained upper bounds 
if we merely rely on upper bound computation.  
Therefore it is desirable to compute 
lower bounds of $L_{2+}$ induced norm efficiently.  
Our upcoming results on such lower bound computation will be reported 
elsewhere in a near future.  

\appendix

\section*{Basics about Copositive Programming}

We first review the definitions and the properties of convex
cones related to the COP.
\begin{definition}\cite{Berman_2003}
The definitions of proper cones 
$\PSD^n$, $\COP^n$, $\CP^n$, $\NN^n$, and $\DNN^n$ in 
$\bbS^n$ are as follows.
\begin{enumerate}
  \item 
  $\PSD^n:=\{P\in\bbS^n:\ \forall x\in\bbR^n,\
	x^TPx\geq 0\}=\{P\in\bbS^n:\ \exists B\
	\mbox{s.t.}\ P=BB^T\}$ is called {\it the positive semidefinite cone}.
  \item
  $\COP^n:=\{P\in\bbS^n:\ \forall x\in\bbR_{+}^n,\
	x^TPx\geq 0\}$ is called {\it the copositive cone}.
  \item
  $\CP^n:=\{P\in\bbS^n:\ \exists B\ge 0\
	\mbox{s.t.}\ P=BB^T\}$ is called {\it the completely positive cone}.
  \item
  $\NN^n:=
  \{P\in\bbS^n:\ P\geq 0\}
  $
  is called {\it the nonnegative cone}.
  \item
  $\PSD^n+\NN^n:=\{ P+Q:\ P\in\PSD^n ,\	Q\in\NN^n \}$.
       This is the Minkowski sum of the positive semidefinite cone and 
       the nonnegative cone.
  \item
  $\DNN^n:=\PSD^n\cap\NN^n$ is called 
  {\it the doubly nonnegative cone}.
\end{enumerate}
\label{def:cone} 
\end{definition}

From Definition \ref{def:cone}, 
we clearly see that the following inclusion relationships hold: 
\begin{equation}
 \CP^n\subset\DNN^n\subset\PSD^n\subset\PSD^n+\NN^n\subset\COP^n,
\label{eq:inc1}
\end{equation}
\begin{equation}
 \CP^n\subset\DNN^n\subset\NN^n\subset\PSD^n+\NN^n\subset\COP^n.
\label{eq:inc2}
\end{equation}
In particular, when $n\leq 4$, 
it is known that  
$\COP^n=\PSD^n+\NN^n$ 
and 
$\CP^n =\DNN^n$  
hold \cite{Berman_2003}. 
On the other hand, as for the duality of these cones, 
$\COP^n$ and $\CP^n$ are dual to each other, 
$\PSD^n+\NN^n$ and $\DNN^n$ are dual to each other, 
and $\PSD^n$ and $\NN^n$ are self-dual.
It is also well known that 
the interiors of the cones $\PSD^n$ and $\NN^n$can be
characterized by
\[
\begin{array}{@{}lcl}
  (\PSD^n)^\circ&=&\{P\in\bbS^n:\  \forall
   x\in\bbR^n\backslash\{0\},\ x^TPx>0\},\\
 (\NN^n)^\circ&=&  \{P\in\bbS^n:\ P>0\}.  
\end{array}
\]

As noted, the COP is a convex optimization problem on the copositive cone, 
and its dual is a convex optimization problem on the completely positive cone.  
As mentioned in \cite{Dur_2010}, 
the problem to determine whether a given symmetric matrix is
copositive or not is a co-NP complete problem, and
the problem to determine whether a given symmetric matrix is
completely positive or not is an NP-hard problem.  
Therefore, it is hard to solve COP 
numerically in general.  
However, since the problem to determine whether a given matrix is in $\PSD+\NN$
or in $\DNN$ can readily be reduced to SDPs, 
we can numerically solve the convex optimization problems on the cones
$\PSD+\NN$ and $\DNN$ easily.  
Moreover, when $n\leq 4$, it is known that $\COP^n=\PSD^n+\NN^n$ and 
$\CP^n=\DNN^n$ as stated above, 
and hence those COPs with $n\leq 4$ can be reduced to SDPs.

\section*{On the Existence of Interior Point Solutions for the SDP
\rec{eq:SDPnew} and Its Dual}

We first note that the dual of the SDP \rec{eq:SDPnew} is given by
\begin{subequations}
\[
\sup_{Z} \ 
\trace\left(
\begin{bmatrix}
C_{a,N}^T \\ D_{a,N}^T
\end{bmatrix}^T Z
\begin{bmatrix}
C_{a,N}^T \\ D_{a,N}^T
\end{bmatrix}\right)\quad \mathrm{subject\ to}
\]
\begin{equation}
Z=
\begin{bmatrix}
Z_{11} & Z_{12} & Z_{13}\\
\ast & Z_{22} & Z_{23}\\
\ast & \ast & Z_{33}\\
\end{bmatrix}
\in \PSD^{n+\np+\nw},
\end{equation}
\begin{equation}
\trace(Z_{33})=1,
\end{equation}
\begin{equation}
A_aZ_a+B_aZ_b^T+(A_aZ_a+B_aZ_b^T)^T=0,
\end{equation}
\begin{equation}
Z_c\ge 0,
\end{equation}
\label{eq:dual}
\end{subequations}
where
\[
\begin{array}{@{}l}
Z_a:=
\begin{bmatrix}
Z_{11} & Z_{12} \\
\ast & Z_{22} 
\end{bmatrix}
\in \PSD^{n+\np},\\   
Z_b:=
\begin{bmatrix}
Z_{13} \\
Z_{23} 
\end{bmatrix}
\in \bbR^{(n+\np)\times \nw},\hspace*{-20mm}\\ 
Z_c:=
\begin{bmatrix}
Z_{22} & Z_{23} \\
\ast & Z_{33} 
\end{bmatrix}
\in \PSD^{\np+\nw}.  
\end{array}
\]

On the existence of the interior point solutions for 
the (primal) SDP \rec{eq:SDPnew} and its dual \rec{eq:dual}, 
we can establish the following two theorems.  

\begin{theorem}
The SDP \rec{eq:SDPnew} has an interior point solution.  
\label{th:primal}
\end{theorem}
\begin{proofof}{\rth{th:primal}}
Let us denote by $P_0\in\bbS_{++}^{n+n_p}$
the unique solution of the Lyapunov equation 
\[
P_0A_a+A_a^TP_0+C_a^TC_a+2 I_{n+\np} =0.  
\]
In addition, with sufficiently small $\varepsilon>0$ we let
\[
 Q_0=I_{\np+\nw}+\varepsilon \one_{\np+\nw}\one_{\np+\nw}^T.  
\]
Then, it is clear that 
\[
 P_0 \in \bbS^{n+\np},\ Q_0 \in (\PSD^{\np+\nw})^\circ+(\NN^{\np+\nw})^\circ.  
\]
Furthermore, if we let $\gamma>0$ sufficiently large, 
we can confirm that  \rec{eq:L2+COPnew} with 
$\preceq$ being replaced by $\prec$ holds with 
$P_a=P_0$ and $Q=Q_0$.  
This clearly shows that the SDP \rec{eq:SDPnew} has an interior point solution.  
\end{proofof}
\begin{theorem}
The dual SDP \rec{eq:dual} has an interior point solution.  
\label{th:dual}
\end{theorem}
We need the next lemma for the proof of \rth{th:dual}.  
\begin{lemma}
For $A_p\in \bbR^{N\nw\times N\nw}$ and $B_p\in \bbR^{N\nw\times \nw}$ 
given by \rec{eq:FIR}, 
the unique solution $\hatZ_p$ to the following Lyapunov equation 
satisfies $\hatZ_p>0$.  
\begin{equation}
A_p \hatZ_p+\hatZ_p A_p^T+ B_p\one_{\nw}\one_{\nw}^T B_p^T=0.  
\label{eq:Lyaplem}
\end{equation}
\label{le:pos}
\end{lemma}
\begin{proofof}{\rle{le:pos}}
From the well-known analytic expression for the solution of Lyapunov equation, 
we have
~\\
\[
\scalebox{0.8}{$
\begin{array}{@{}lcl}
\hatZ_p&=&\int_0^\infty \exp(A_pt) B_p\one_{\nw}\one_{\nw}^T B_p^T \exp(A_p^T t) dt \\
&=&
\displaystyle\int_0^\infty 
\begin{bmatrix}
\frac{t^{N-1}}{(N-1)!} \exp(\alpha t)  I_{\nw} \\ \vdots \\ \exp(\alpha t)  I_{\nw}
\end{bmatrix}
\one_{\nw}\one_{\nw}^T \\
& &
\begin{bmatrix}
\frac{t^{N-1}}{(N-1)!} \exp(\alpha t)  I_{\nw} \\ \vdots \\ \exp(\alpha t)  I_{\nw}
\end{bmatrix}^T dt.  
\end{array}$}
\]
~\\~\\
This clearly shows that $\hatZ_p>0$.  
\end{proofof}

We are now ready to prove \rth{th:dual}. 

\begin{proofof}{\rth{th:dual}}
We first denote by $Z_{a,0}\in \bbS_+^{n+\np}$
the unique solution to the Lyapunov equation given below
that is obtained by substituting $Z_b=B_a$ in (\ref{eq:dual}c):  
\[
A_a Z_{a,0}+Z_{a,0}A_a^T+2B_aB_a^T=0.  
\]
Then, since the pair $(A_a,B_a)$ is controllable, 
we see $Z_{a,0}\succ 0$.  
With this fact in mind, we next substitute 
\begin{equation}
Z_b=B_a+
\begin{bmatrix}  0_{n,\nw} \\ \varepsilon \one_{\np}\one_{\nw}^T \end{bmatrix}
=:Z_{b,\varepsilon}
\left(
=: 
\begin{bmatrix}  \hatZ_{13} \\ \hatZ_{23} \end{bmatrix}
\right)
\label{eq:Zb}
\end{equation}
in (\ref{eq:dual}c) with sufficiently small $\varepsilon>0$
and consider the resulting Lyapunov equation:  
\begin{equation}
A_a Z_{a,\varepsilon}+Z_{a,\varepsilon}A_a^T+B_{a}Z_{b,\varepsilon}^T+Z_{b,\varepsilon}B_{a}^T=0.  
\label{eq:Lyape}
\end{equation}
Then, for the unique solution  $Z_{a,\varepsilon}\in\bbS^{n+\np}$
to this Lyapunov equation, 
we see that $Z_{a,\varepsilon}\succ 0$ holds for sufficiently small $\varepsilon>0$
since $Z_{a,0}\succ 0$.  
If we partition $Z_{a,\varepsilon}\succ 0$ as
\begin{equation}
Z_{a,\varepsilon}=:
\begin{bmatrix}
\hatZ_{11}  & \hatZ_{12}  \\ \ast  & \hatZ_{22}  \\
\end{bmatrix}\succ 0,\ \hatZ_{22}\in \PSD^{\np},  
\label{eq:Za}
\end{equation}
we see that the $(2,2)$-block of \rec{eq:Lyape} can be written as
\[
\begin{array}{@{}l}
A_p \hatZ_{22} +\hatZ_{22}  A_p\\
+B_p(B_p+\varepsilon \one_{\np}\one_{\nw}^T)^T+(B_p+\varepsilon \one_{\np}\one_{\nw}^T)B_p^T=0.  
\end{array}
\]
If we compare the above equation with \rec{eq:Lyaplem} in \rle{le:pos}, 
we see 
\[
\begin{array}{@{}l}
B_p(B_p+\varepsilon \one_{\np}\one_{\nw}^T)^T+
(B_p+\varepsilon \one_{\np}\one_{\nw}^T)B_p^T \\
\ge \varepsilon (B_p \one_{\nw}\one_{\np}^T +\one_{\np}\one_{\nw}^T B_p^T )\\
\ge \varepsilon (B_p \one_{\nw}\one_{\nw}^T B_p^T +B_p\one_{\nw}\one_{\nw}^T B_p^T )\\
=2 \varepsilon B_p \one_{\nw}\one_{\nw}^T B_p^T.  
\end{array}
\]
It follows that 
\begin{equation}
\hatZ_{22}\ge 2\varepsilon \hatZ_p>0.  
\label{eq:Z22} 
\end{equation}
Finally, if we let
\begin{equation}
Z_{\nu,33}:=\nu I_{\nw}+\one_{\nw}\one_{\nw}^T (=\hatZ_{33})
\end{equation}
with sufficiently large $\nu>0$, we see from \rec{eq:Za} that
\[
 \hatZ=
\begin{bmatrix}
\hatZ_{11} & \hatZ_{12} & \hatZ_{13}\\
\ast & \hatZ_{22} & \hatZ_{23}\\
\ast & \ast & \hatZ_{33}\\
\end{bmatrix}\succ 0.  
\]
In addition, if we define
\[
Z:=\frac{1}{\trace(\hatZ_{33})} \hatZ, 
\]
we can conclude that this $Z$ satisfies 
$Z\succ 0$, (\ref{eq:dual}b), (\ref{eq:dual}c), and $Z_c>0$.  
To summarize, the dual SDP \rec{eq:dual} has an interior point solution.  
\end{proofof}
%


\end{document}